\documentclass[12pt]{article}
\makeatletter
\newfont{\footsc}{cmcsc10 at 8truept}
\newfont{\footbf}{cmbx10 at 8truept}
\newfont{\footrm}{cmr10 at 10truept}
\makeatother
\usepackage{amssymb, amsmath, amsthm}
\usepackage{epsfig}

\newtheorem{theorem}{Theorem}
\newtheorem{prop}{Proposition}
\newtheorem{lem}{Lemma}

\newcommand{\Bild}[4]{
\begin{figure}[tb]
  \begin{center}
    \leavevmode
    \epsfig{file=#2,height=#1cm}\vspace{-0.4cm}
    \caption{{\small #3}}
    \label{#4}
  \end{center}
\end{figure}}

\title{{\bf On the Entropy and Letter Frequencies\\ 
of Ternary Square-Free Words}}

\author{
\begin{tabular}{cc}
{\bf Christoph~Richard} & 
  {\bf Uwe~Grimm} \\
{\small Institut f\"ur Mathematik}&
  {\small Applied Mathematics Department}\\[-0.8ex]
{\small Universit\"at Greifswald}&
  {\small The Open University}\\[-0.8ex]
{\small Jahnstr. 15a}&
  {\small Walton Hall}\\[-0.8ex]
{\small 17487 Greifswald, Germany}& 
  {\small Milton Keynes MK7 6AA, UK}\\[-0.8ex]
{\small \texttt{richard@uni-greifswald.de}} &
  {\small \texttt{u.g.grimm@open.ac.uk}}\\
\end{tabular}
}

\date{\small \today\\ %Submitted: \today;  Accepted: \today.}\\
}

\begin{document}
\maketitle

\begin{abstract}
We enumerate all ternary length-$\ell$ square-free words, which are
words avoiding squares of words up to length $\ell$, for $\ell\le 24$.
We analyse the singular behaviour of the corresponding generating
functions. This leads to new upper entropy bounds for ternary
square-free words.  We then consider ternary square-free words with
fixed letter densities, thereby proving exponential growth for certain
ensembles with various letter densities.  We derive consequences for
the free energy and entropy of ternary square-free words.
\end{abstract}

\section{Introduction}

The interest in the combinatorics of pattern-avoiding
\cite{BEN79,BNT89,C93}, in particular of power-free words, goes back
to work of Axel Thue in the early 20th century \cite{Thue1,Thue2}. The
celebrated Prouhet-Thue-Morse sequence, defined by a substitution rule
$a\rightarrow ab$ and $b\rightarrow ba$ on a two-letter alphabet
$\{a,b\}$, proves the existence of infinite cube-free words in two
letters $a$ and $b$.

Here, a word of length $n$ is a string of $n$ letters from a certain
alphabet $\Sigma$, an element of the language
$\mathcal{L}(n)=\Sigma^n$ of $n$-letter words in $\Sigma$. The union
\begin{equation}
\mathcal{L}\;=\;\bigcup_{n\ge 0}\mathcal{L}(n)
\;=\;\Sigma^{\mathbb{N}_0}
\end{equation}
is the language of all words in the alphabet $\Sigma$. It is a monoid,
with concatenation of words as operation, and with the empty word
$\lambda$ of zero length as neutral element \cite{L83}. A word $w$ is
called {\it square-free}\/ if $w=xyyz$, with words $x$, $y$ and $z$,
implies that $y=\lambda$ is the empty word, and cube-free words are
defined analogously. So square-free words are characterised by the
property that they do not contain an adjacent repetition of any
subword.

It is easy to see that there are only a few square-free words in two
letters, these are the empty word $\lambda$, the two letters $a$ and
$b$, the two-letter words $ab$ and $ba$, and, finally, the
three-letter words $aba$ and $bab$. Appending any letter to those two
words inevitably results in a square, either of a single letter, or of
one of the square-free two-letter words.

However, there do exist infinite ternary square-free words, i.e.,
square-free words on a three-letter alphabet. In fact, the number
$s_{n}$ of ternary square-free words of length $n$ grows exponentially
with $n$. Denoting set sets of ternary square-free words of length $n$
by $\mathcal{A}_{n}$, we have
\begin{eqnarray}
\mathcal{A}_0 & = & \{\lambda\},\nonumber\\
\mathcal{A}_1 & = & \{a,b,c\},\nonumber\\
\mathcal{A}_2 & = & \{ab,ac,ba,bc,ca,cb\},\nonumber\\
\mathcal{A}_3 & = & \{aba,abc,aca,acb,bab,bac,bca,bcb,cab,cac,cba,cbc\},
\end{eqnarray}
and so on. So $s_0=1$, $s_1=3$, $s_2=6$, $s_3=12$, and so on, see
\cite{BEG97} and \cite{G01} where the values of $s_{n}$ for $n\le 90$
and $91\le n\le 110$ are tabulated, respectively. In \cite{SP}, the
sequence $s_{n}$ is listed as A006156 (formerly M2550).

In this article, we consider ternary square-free words
\cite{Thue1,Thue2,Z58,P70,BEN79,Bra83,Bri83,E83,L83,S83,L85,S85,K86,BEG97,%
EZ98,NZ99,G01,C02,S02}. We are interested in the asymptotic growth of
the sequence $s_{n}$. We use a series of generating functions for a
truncated square-freeness condition and conjecture the presence of a
natural boundary at the radius of convergence. We also consider the
frequencies of letters in ternary square-free words and derive upper
and lower bounds.  We prove exponential growth for certain ensembles
of ternary square-free words with fixed letter frequencies.  We use
methods of statistical mechanics \cite{J00} to prove that, subject to
a plausible regularity assumption on the free energy of ternary
square-free words, the maximal exponential growth occurs for words
with equal mean letter frequencies, where we average over all
square-free words.  Some of our results are based on extensive exact
enumerations of square-free ternary words of length $n\le 110$
\cite{G01} and on constructions of generalised Brinkhuis triples
\cite{E83,G01}.

\section{Ternary square-free words}

Denote the number of ternary square-free words by $s_n$ and the
corresponding generating function by $S(x)$,
\begin{equation}
S(x)\;=\; \sum_{n=0}^{\infty}\, s_{n}\, x^{n}\,.
\end{equation}
Since the language of ternary square-free words is subword-closed, we
conclude that the sequence $s_n$ is submultiplicative,
\begin{equation}
s_{n+m} \;\le\; s_n \, s_m\,.
\end{equation}
A standard argument, compare \cite[Lemma 1]{BEG97} and \cite[Lemma
A.1]{J00}, shows that this guarantees that the limit
$\mathcal{S}:=\lim_{n \to \infty} \frac{1}{n} \log s_n$, also called
the {\it entropy}, exists.  Bounds for the limit have been obtained in
a number of investigations \cite{Bri83,Bra83,E83,EZ98,NZ99,G01,S02},
which give
\begin{equation}
1.1184\;\approx\; 110^{1/42}\;\le\; \exp(\mathcal{S})\;<\;1.30201064\,,
\label{bounds}
\end{equation}
but the exact value is unknown.  The lower bound implies an
exponential growth of $s_n$ with $n$. The behaviour of the subleading
corrections to the exponential growth is not understood.

One of the authors computed the numbers $s_n$ for $n\le 110$
\cite{G01}. Assuming an asymptotic growth of the numbers $s_n$ of the
form
\begin{equation}\label{eqn:asympt}
s_n\;\sim\; A\, x_{c}^{-n}\, n^{\gamma-1} \qquad (n \to \infty)\,,
\end{equation}
we used differential approximants \cite{G89} of first order to get
estimates of the critical point $x_c=\exp(-\mathcal{S})$, the critical
exponent $\gamma$ and the critical amplitude $A$.  We obtain
\begin{equation}
A \;=\; 12.72(1)\,, \qquad 
x_c\; =\; 0.768189(1)\,, \qquad 
\gamma \;=\; 1.0000(1)\,,
\end{equation}
where the number in the bracket denotes the (estimated) uncertainty in
the last digit. The value of $\gamma$, also found in \cite{NZ99},
suggests a simple pole as dominant singularity of the generating
function at $x=x_c$. Numerical analysis indicates the presence of a
natural boundary, a topic which we considered further by computing
approximating generating functions $S^{(\ell)}(x)$, which count the
number of words which contain no squares of words of length $\le\ell$.

\section{Generating functions}

We call a word $w\in{\cal L}$ {\it length-$\ell$ square-free} if
$w=xyyz$, with $x,z\in{\cal L}$ and $y\in\bigcup_{n=0}^{\ell}{\cal
L}(n)$, implies that $y$ is the empty word $\lambda$. In other words,
$w$ does not contain the square of a word of length $\le\ell$.

Denote the number of ternary length-$\ell$ square-free words of length
$n$ by $s_{n}^{(\ell)}$. Clearly, $\ell'>\ell$ implies
$s_{n}^{(\ell')}\le s_{n}^{(\ell)}$, because at least the same number
of words are excluded. On the other hand, we have
$s_{n}^{(\ell')}=s_{n}^{(\ell)}=s_{n}$ for $n<2\ell$.  Thus, by
considering larger and larger squares $\ell$, we approach the case of
square-free words.

We define corresponding generating functions
\begin{equation}
S^{(\ell)}(x) \;=\; \sum_{n=0}^{\infty}\,s_{n}^{(\ell)}\, x^{n}
\end{equation}
for the number of ternary length-$\ell$ square-free words. These
generating functions are rational functions of the variable $x$ which
can be calculated explicitly, at least for small values of $\ell$, see
\cite{NZ99} where the computation is explained in detail.  The first
few generating functions are
\begin{eqnarray*}
S^{(0)}(x) &\!=\!& \frac{1}{1-3x}\,, \\
S^{(1)}(x) &\!=\!& \frac{1+x}{1-2x}\,, \\
S^{(2)}(x) &\!=\!& \frac{1+2x+2x^{2}+3x^{3}}{1-x-x^{2}}\,, \\
S^{(3)}(x) &\!=\!& \frac{1\!+\!3x\!+\!6x^{2}\!+\!11x^{3}\!+\!14x^{4}\!+\!
                     20x^{5}\!+\!20x^{6}\!+\!21x^{7}\!+\!12x^{8}\!+\!
                     6x^{9}(1\!-\!x\!-\!x^{2}\!-\!x^{3}\!-\!x^{4})}
                    {1-x^{3}-x^{4}-x^{5}-x^{6}}\,.
\end{eqnarray*}
We computed the generating functions $S^{(\ell)}(x)$ explicitly for
$\ell\le 24$. The functions are available as Mathematica code \cite{W}
at \cite{G03}. Note that some generating functions agree; for
instance, $S^{(4)}(x)=S^{(5)}(x)$. The reason is that, going from
$\ell=4$ to $\ell=5$, no ``new'' squares arise; in other words, all
squares of square-free words of length $5$ already contain a square of
a word of smaller length.

The radius of convergence $x_{c}^{(\ell)}\le x_{c}$ of the series
defining the generating function $S^{(\ell)}(x)$ is determined by a
pole in the complex plane located closest to the origin, thus by a
zero of the denominator polynomial of smallest modulus.  Due to
Pringsheim's theorem \cite[Sec.~7.21]{St78}, a real and positive such
zero exists.  Note that the zeros of the numerator and denominator are
mutually exclusive, because the do not contain common polynomial
factors.

\begin{table}[tb]
\caption{Degrees $d_{\rm num}$ and $d_{\rm den}$ of the numerator and
denominator polynomials of the generating functions $S^{(\ell)}(x)$,
respectively, and the numerical values of the radius of convergence
$x_{c}^{(\ell)}$.\label{tab1}}\smallskip
\begin{center}
\begin{tabular}{r@{\qquad}r@{\qquad}r@{\qquad}l}
\hline
\multicolumn{1}{c}{$\ell$\rule[-1ex]{0ex}{4ex}} & 
\multicolumn{1}{c}{$d_{\rm num}$} & 
\multicolumn{1}{c}{$d_{\rm den}$} & 
\multicolumn{1}{c}{$x_{c}^{(\ell)}$} \\
\hline\rule[0ex]{0ex}{3ex}
 $0$            &    $0$ &    $1$ & $0.333\,333\,333$ \\
 $1$            &    $1$ &    $1$ & $0.500\,000\,000$ \\
 $2$            &    $3$ &    $2$ & $0.618\,033\,989$ \\
 $3$            &    $5$ &    $3$ & $0.682\,327\,804$ \\
 $4$, $5$       &   $13$ &    $6$ & $0.724\,491\,959$ \\
 $6$, $7$       &   $27$ &   $15$ & $0.750\,653\,202$ \\
 $8$, $9$, $10$ &   $38$ &   $19$ & $0.757\,826\,433$ \\
$11$            &   $81$ &   $58$ & $0.762\,463\,266$ \\
$12$            &  $143$ &  $106$ & $0.765\,262\,611$ \\
$13$, $14$      &  $184$ &  $145$ & $0.766\,784\,948$ \\
$15$            &  $209$ &  $170$ & $0.767\,006\,554$ \\
$16$, $17$      &  $217$ &  $178$ & $0.767\,136\,379$ \\
$18$            &  $441$ &  $380$ & $0.767\,542\,044$ \\
$19$            &  $644$ &  $594$ & $0.767\,752\,831$ \\
$20$            &  $968$ &  $890$ & $0.767\,887\,486$ \\
$21$            & $1003$ &  $925$ & $0.767\,896\,727$ \\
$22$            & $1436$ & $1337$ & $0.767\,974\,175$ \\
$23$            & $1966$ & $1872$ & $0.768\,042\,881$ \\
\rule[-1ex]{0ex}{1ex}
$24$            & $2905$ & $2787$ & $0.768\,085\,659$ \\
\hline
\end{tabular}
\end{center}
\end{table}

The values $x_c^{(\ell)}$ are given in Table~\ref{tab1}, together with
the degrees $d_{\rm num}$ and $d_{\rm den}$ of the polynomials in the
numerator and in the denominator which both grow with $\ell$. Thus,
with growing length $\ell$, the generating functions $S^{(\ell)}(x)$
have an increasing number of zeros and poles. The patterns of zeros
and poles appear to accumulate in the complex plane close to the unit
circle around the origin; and comparing the patterns for increasing
$\ell$ one might be tempted to the plausible conjecture that the poles
approach the unit circle in the limit as
$\ell\rightarrow\infty$. However, there appear to be some oscillations
in the patterns close to the real line, and at present we dot not have
any argument why the poles should accumulate on the unit circle.

The values $x_{c}^{(\ell)}$ in Table~\ref{tab1} approach $x_{c}$ from
below, so they yield upper bounds on the exponential growth constant
$\mathcal{S}=-\log(x_{c})$. The upper bound quoted in equation
(\ref{bounds}) above was given in \cite{NZ99} on the basis of an
estimate for $x_{c}^{(23)}$ obtained via the series expansion of
$S^{(23)}(x)$. Our value for $x_{c}^{(23)}$, based on the complete
evaluation of the generating function $S^{(23)}(x)$, is contained in
Table~\ref{tab1}; it confirms the bound of Noonan and Zeilberger
\cite{NZ99}. The value for $\ell=24$ slightly improves the upper
bound.
\begin{theorem}\label{th:uppbound}
The entropy $\mathcal{S}$ of ternary square-free words is bounded as\/
$\mathcal{S}\le -\log(x_{c}^{(24)})$, which gives\/
$\exp(\mathcal{S})<1.30193812<1/x_{c}^{(24)}$.
\qed 
\end{theorem}

The complete set of poles of the generating function $S^{(24)}(x)$ is
shown in Fig.~\ref{fig:poles}. The pattern looks very similar for
other values of $\ell$. This suggests that, in the limit as $\ell$
becomes infinite, which corresponds to the generating function $S(x)$
of ternary square-free words, the poles accumulate close to the unit
circle. This corroborates the conjecture that $S(x)$ has a natural
boundary.

\Bild{9}{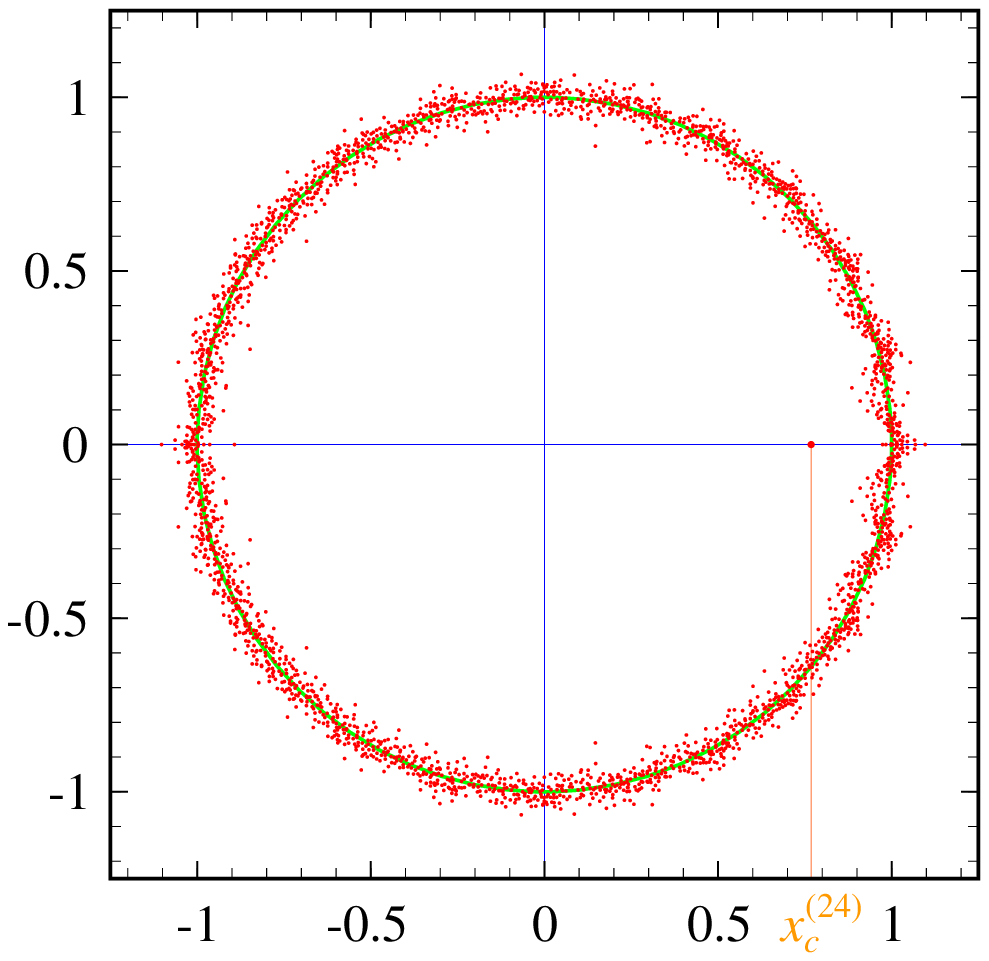}{Pattern of poles of the generating functions
              $S^{(24)}(x)$ in the complex plane. The poles (red)
              accumulate along the unit circle (green). The isolated
              pole at $x^{(24)}_c$ on the real positive axis determines
              the radius of convergence.}{fig:poles}

\section{Square-free words with fixed letter frequencies}

We now consider the letter statistics of ternary square-free words.
Denote the number of occurrences of the letter $a$ in a ternary
square-free word $w_{n}$ of finite length $n$ by $a(w_n)$. Clearly,
the frequency of the letter $a$ in $w_{n}$ is $0\le a(w_{n})/n\le 1$.
For an infinite ternary square-free word $w$, letter frequencies do
not generally exist. Consider sequences $\{w_{n}\}$ of $n$-letter
subwords containing arbitrarily long words. We define upper and lower
frequencies $f_a^{+}\ge f_a^{-}$ by
$f_a^{+}:=\sup_{\{w_{n}\}}\limsup_{n\to\infty}a(w_n)/n$ and
$f_a^{-}:=\inf_{\{w_{n}\}}\liminf_{n\to\infty}a(w_n)/n$, where we take
the supremum and infimum over all sequences $\{w_{n}\}$.  We can also
compute these from $a_n^{+}=\max_{w_n\subset w} a(w_n)$ and
$a_n^{-}=\min_{w_n\subset w} a(w_n)$ by
$f_a^{\pm}=\lim_{n\to\infty}a_{n}^{\pm}/n$, as these limits exist.
This follows, for instance, from the subadditivity of the sequences
$\{a_n^{+}\}$ and $\{1-a_{n}^{-}\}$.  If the infinite word $w$ is such
that $f_a^{+}=f_{a}^{-}=:f_{a}$, we call $f_{a}$ the {\it frequency}\/
of the letter $a$ in $w$. In general, $f_a^{+}>f_a^{-}$, and letter
frequencies do not exist, see also the discussion below.

However, we can derive bounds on the upper and lower letter
frequencies $f_a^{+}$ and $f_a^{-}$.  Denote the number of ternary
square-free words of length $n$ which contain the letter $a$ exactly
$k$ times by $s_{n,k}$.  Since there are no square-free words of
length greater than three in two letters, a ternary square-free word
contains no gaps between letters $a$ of length greater than three.
This implies $s_{n,k}=0$ for $k<n/4$ or $k>n/2$, since the minimal
number of letters $b$ and $c$ is, by the same argument, equal to
$k=n/2$.  By counting the number $s_{n,k}$ of ternary square-free
words with a given number $k$ of letters $a$, we can sharpen these
bounds. Clearly, for fixed $k$, there are numbers $n_{\rm min}(k)$ and
$n_{\rm max}(k)$ such that $s_{n,k}=0$ for $n<n_{\rm min}(k)$ and
$n>n_{\rm max}(k)$. This means that any ternary square-free word of
length $(m+1)n_{\rm max}(k)\ge n>m n_{\rm max}(k)$, for any integer
$m$, contains at least $m k+1$ letters $a$, so the frequency of the
letter $a$ is bounded from below by $(m k+1)/(m n_{\rm max}(k)+1)$,
which becomes $k/n_{\rm max}(k)$ as $m$ tends to infinity. Similarly,
any word of length $m n_{\rm min}(k)> n\ge (m-1)n_{\rm min}(k)$
contains at most $mk-1$ letters $a$. Thus we obtain an upper limit of
$(m k-1)/(m n_{\rm min}(k)-1)$, which becomes $k/n_{\rm min}(k)$ as
$m$ tends to infinity. We computed $n_{\rm max}(k)$ for $k\le 31$ and
$n_{\rm min}(k)$ for $k\le 40$; the strongest bounds are derived from
$n_{\rm max}(31)=117$ and $n_{\rm min}(39)=97$, which yield lower and
upper bounds $31/117\approx 0.265$ and $39/97\approx 0.402$,
respectively, for the frequency of a single letter in an infinite
ternary square-free word.  This gives
\begin{theorem}\label{th:fbounds}
The upper and lower frequencies\/ $f^{\pm}$ of a given letter in an
infinite ternary square-free word are bounded by\/ $0.265\approx31/117
\le f^{-}\le f^{+}\le 39/97\approx 0.402$. \qed
\end{theorem}
\noindent {\bf Remark.} In fact, there is a recent, stronger result
for the lower frequency \cite{Tara02}.  The minimum frequency
$f^{-}_{\rm min}$ is bounded from below and above by \cite{Tara02}
\[
0.274649\approx 1780/6481\le f^{-}_{\rm min}\le 64/233\approx 
0.274678\, ,
\]
compare also similar treatments for binary power-free words
\cite{KK97,KKT98}. The upper bound can be sharpened to $f^{+}\le
469/1201\approx 0.390508$ \cite{Tara03}.
\vspace{2ex}

It is easy to see that the {\it mean letter frequency}\/ of any given
letter in the set of ternary square-free words is $1/3$. This is a
consequence of symmetry under permutation of letters.  Indeed, the
symmetric group $S_{3}$ acts on any square-free word $w$ by
permutation of the three letters, and the set of square-free words of
a given length is a disjoint union of orbits under this action. Each
orbit consists of a square-free word and its images under permutation
of letters, and each letter has the same mean frequency on this
orbit. So, for each orbit, the mean frequency of any given letter is
$1/3$, thus also for the set of all ternary square free words of any
given length, or indeed for the set of all ternary square free words.

We now want to show that there exist ternary square-free words of
infinite length with well-defined letter frequencies for the case
$f_a=f_b=f_c=1/3$ and for some cases where not all letter are equally
frequent. In fact, we are going to prove not just that, but that there
are exponentially many such words, so the growth rate for words of
fixed frequencies, at least for the cases considered below, is
positive. This can be done in a similar fashion as the proofs that the
number of ternary square-free words grow exponentially
\cite{Bri83,Bra83,E83,EZ98,NZ99,G01,S02}. These proofs are based on
Brinkhuis triple pairs \cite{Bri83,Bra83,E83,EZ98,NZ99} and their
generalisations \cite{E83,G01,S02}. We briefly sketch the argument
here, see \cite{Bri83,Bra83,E83,EZ98,NZ99,G01,S02} for details.

The argument is based on square-free morphisms \cite{C82,C83}.  Here,
we immediately consider the generalised version of
\cite{E83,G01}. Assume that we have a set of substitution rules
\begin{equation}
a \;\rightarrow\; \begin{cases}
w_{a}^{(1)}\\
w_{a}^{(2)}\\
\vdots\\
w_{a}^{(k)}
\end{cases}\qquad
b \;\rightarrow\; \begin{cases}
w_{b}^{(1)}\\
w_{b}^{(2)}\\
\vdots\\
w_{b}^{(k)}
\end{cases}\qquad
c \;\rightarrow\; \begin{cases}
w_{c}^{(1)}\\
w_{c}^{(2)}\\
\vdots\\
w_{c}^{(k)}
\end{cases}
\label{brinktrip}
\end{equation}
where $w_{a}^{(j)}$, $w_{b}^{(j)}$ and $w_{c}^{(j)}$, $1\le j\le k$,
are ternary square-free words of equal length $m$. Starting from any
ternary square-free word $w$ of length $n$, consider the set of all
words of length $mn$ obtained by substituting each letter, choosing
independently one of the $k$ words from the lists above. A {\it
generalised Brinkhuis triple}\/ is defined as a set of substitution
rules (\ref{brinktrip}) such that all these words of length $mn$ are
square-free, for any choice of $w$. This immediately implies that the
number of square-free words grows at least as $k^{1/(m-1)}$, see
\cite[Lemma~2]{G01}. In the case $k=1$, this reduces to a usual
substitution rule without any freedom; in this case, it only proves
existence of infinite words, not exponential growth of the number of
words with length.

In \cite{G01}, a special class of generalised Brinkhuis triples was
considered, and triples up to length $m=41$ with $k=65$ were
obtained. This was recently improved to $m=43$ and $k=110$ in
\cite{S02}, yielding the lower bound of (\ref{bounds}).

What about the letter frequencies? In general, the words $w_{a}^{(j)}$
that replace $a$ will have different letter frequencies, and in this
case it is easy to see that not all the infinite words obtained by
repeated substitution will have well-defined letter
frequencies. However, we can say something about letter frequencies if
we consider generalised Brinkhuis triples where all words
$w_{a}^{(j)}$, $1\le j\le k$, have the {\em same}\/ letter
frequencies, and analogously for the words $w_{b}^{(j)}$, $1\le j\le
k$, and $w_{c}^{(j)}$, $1\le j\le k$. In this case, regardless of our
choice of words in the substitution process, we obtain words with
well-defined letter frequencies, precisely as in the case of a
standard substitution rule. Denoting the number of letters $a$, $b$
and $c$ in any of the words $w_{a}^{(j)}$ by $n_a^a$, $n_a^b$ and
$n_a^c$, respectively, with $n_a^a+n_a^b+n_a^c=m$, and analogously for
$w_{b}^{(j)}$ and $w_{c}^{(j)}$, we can summarise the letter-counting
for the generalised Brinkhuis triple in a $3\times 3$ substitution
matrix
\begin{equation}
M \;=\; \begin{pmatrix}
n_a^a & n_b^a & n_c^a\\
n_a^b & n_b^b & n_c^b\\
n_a^c & n_b^c & n_c^c
\end{pmatrix}\, .
\end{equation}
In general, all entries of this matrix are positive integers, because
there are no square-free words of length $m>3$ with only two letters.
The (right) Perron-Frobenius eigenvector is thus positive, and its
components encode the letter frequencies of the infinite words
obtained by repeated application of the substitution rules. The
Perron-Frobenius eigenvalue is $m$, because $(1,1,1)$ is a
left eigenvector with eigenvalue $m$.

As mentioned previously, the generalised Brinkhuis triples considered
in \cite{G01} do not have the property that the letter frequencies of
the substitution words coincide. However, if we have a generalised
Brinkhuis triple, any subset of substitutions also forms a triple,
because all we do is restricting to a subset of words which still are
square-free. So by looking at the triples of \cite{G01} and selecting
suitable subsets of substitutions, we can use the same arguments to
prove exponential growth of words with fixed letter frequencies.

\subsection{Equal letter frequencies}

Let us first consider the case of equal frequencies
$f_{a}=f_{b}=f_{c}=1/3$. We note that the special Brinkhuis triples of
\cite{G01} had the additional property that
$w_b^{(j)}=\sigma(w_a^{(j)})$ and $w_c^{(j)}=\sigma^{2}(w_a^{(j)})$,
where $\sigma$ is the permutation of letters defined by $\sigma(a)=b$
and $\sigma(b)=c$. If we select a subset of the words replacing $a$
such that they have the same numbers of letters $n_a^a$, $n_a^b$ and
$n_a^c$, the substitution matrix for the corresponding triple
consisting of those words and their images under $\sigma$ is
\begin{equation}
M \;=\; \begin{pmatrix}
n_a^a & n_a^c & n_a^b\\
n_a^b & n_a^a & n_a^c\\
n_a^c & n_a^b & n_a^a
\end{pmatrix}
\end{equation}
which is symmetric. Hence the right Perron-Frobenius eigenvector is
$(1,1,1)^{t}$, and the letter frequencies are given by
$f_{a}=f_{b}=f_{c}=1/3$.

The simplest example is a Brinkhuis triple with $m=18$ \cite{G01} (see
also \cite{NZ99}) which explicitly given by
\begin{equation}
\begin{split}
w_{a}^{(1)}&\;=\;abcacbacabacbcacba\,,\\
w_{a}^{(2)}&\;=\;abcacbcabacabcacba\;=\;\overline{w}_{a}^{(1)}\,, 
\end{split}
\end{equation}
where $\overline{w}_{a}^{(1)}$ denotes $w_{a}^{(1)}$ read
back-to-front, which thus has the same letter numbers $n_a^a=7$,
$n_a^b=5$ and $n_a^c=6$. So the number of ternary square-free words
with letter frequencies $f_{a}=f_{b}=f_{c}=1/3$ grows at least as
$2^{1/17}$. By looking for the largest subsets of words with equal
letter frequencies in the special Brinkhuis triples of \cite{G01}, we
can improve this bound. For $m=41$, we find $30$ words $w_{a}^{(j)}$
with letter numbers $n_a^a=14$, $n_a^b=13$ and $n_a^c=14$, yielding a
lower bound of $30^{1/40}\approx 1.08875$ for the exponential of the
entropy. One of the two triples for $m=43$ of \cite{S02} contains $39$
words with $n_a^a=14$, $n_a^b=14$ and $n_a^c=15$.  This gives the
following result.
\begin{lem}\label{th:lowbound1}
The entropy\/ $\mathcal{S}(\tfrac{1}{3},\tfrac{1}{3},\tfrac{1}{3})$ of
ternary square-free words with letter frequencies\/
$f_{a}=f_{b}=f_{c}=1/3$ is bounded from below via\/
$\exp[\mathcal{S}(\tfrac{1}{3},\tfrac{1}{3},\tfrac{1}{3})]\ge
39^{1/42}\approx 1.09115$.  \qed
\end{lem}
\noindent {\bf Remark.} This bound can without doubt be improved,
because the triples of \cite{G01} and \cite{S02} where not optimised
to contain the largest number of words of equal frequency.
\vspace{2ex}

\subsection{Unequal letter frequencies}

What about words with non-equal letter frequencies?  The following
square-free substitution rule \cite{Z58}
\begin{equation}\label{form:subs}
\begin{split}
a &\;\rightarrow\; cacbcabacbab \\
b &\;\rightarrow\; cabacbcacbab \\
c &\;\rightarrow\; cbacbcabcbab
\end{split}
\end{equation}
already shows that infinite words with unequal letter frequencies
exist.  In this case, the substitution matrix is
\begin{equation}
M=\begin{pmatrix}
4 & 4 & 3\\
4 & 4 & 5\\
4 & 4 & 4
\end{pmatrix}\, ,
\end{equation}
and the right Perron-Frobenius eigenvector with eigenvalues $12$ is
$(11,13,12)^{t}$. Thus this substitution leads to a ternary
square-free word with letter frequencies $f_a=11/36$, $f_b=13/36$ and
$f_c=1/3$.

Can we show that, for some frequencies, there are exponentially many
words?  Indeed, for some examples we can find generalised Brinkhuis
triples by choosing subsets of those given in \cite{G01}.  Here, we
restrict ourselves to a few examples.

Consider the two generating words
\begin{equation}
\begin{split}
w_{1}&\;=\;abcbacabacbcabacabcbacbcabcba\qquad\hphantom{c}
\mbox{($n_{a}=10$, $n_{b}=10$, $n_{c}=9$)}\,,\\
w_{2}&\;=\;abcbacabacbcacbacabcacbcabcba\qquad\hphantom{b}
\mbox{($n_{a}=10$, $n_{b}=9$, $n_{c}=10$)} \,,
\end{split}
\end{equation}
of a Brinkhuis triple with $m=29$ \cite{G01}.  Choosing
$w_{a}^{(1)}=w_{1}$, $w_{a}^{(2)}=\overline{w}_{1}$,
$w_{b}^{(1)}=\sigma(w_1)$, $w_{b}^{(2)}=\sigma(\overline{w}_1)$,
$w_{c}^{(1)}=\sigma^{2}(w_{2})$ and
$w_{c}^{(2)}=\sigma^{2}(\overline{w}_{2})$, where again $\overline{w}$
denotes the words obtained by reversing $w$, we obtain a Brinkhuis
triple with substitution matrix
\begin{equation}
M\;=\;
\begin{pmatrix}
10 & 9  & 9 \\  
10 & 10 & 10 \\
9  & 10 & 10
\end{pmatrix}\, .
\end{equation}
The corresponding frequencies are
$f=(f_{a},f_{b},f_{c})=(\tfrac{9}{28},\tfrac{10}{29},\tfrac{271}{812})$,
and the growth rate for this case is at least $2^{1/28}$.

Consider now two generating words
\begin{equation}
\begin{split}
w_{1}&\;=\;abcbacabacbabcabacabcacbcabcba\qquad\hphantom{c}
\mbox{($n_{a}=11$, $n_{b}=10$, $n_{c}=9$)}\,,\\
w_{2}&\;=\;abcbacabacbcabcbacabcacbcabcba\qquad\hphantom{a}
\mbox{($n_{a}=10$, $n_{b}=10$, $n_{c}=10$)}\,,
\end{split}
\end{equation}
of a Brinkhuis triple with $m=30$ \cite{G01}.  Choosing 
$w_{a}^{(1)}=w_{1}$,
$w_{a}^{(2)}=\overline{w}_{1}$, 
$w_{b}^{(1)}=\sigma(w_2)$,
$w_{b}^{(2)}=\sigma(\overline{w}_2)$,
$w_{c}^{(1)}=\sigma^{2}(w_{\alpha})$ and
$w_{c}^{(2)}=\sigma^{2}(\overline{w}_{\alpha})$, 
where $\alpha\in\{1,2\}$, we obtain two Brinkhuis
triples with substitution matrices $M_{\alpha}$ given by
\begin{equation}
M_{1}\;=\;
\begin{pmatrix}
11 & 10 & 10 \\  
10 & 10 & 9  \\
9  & 10 & 11
\end{pmatrix}\, ,\qquad
M_{2}\;=\;
\begin{pmatrix}
11 & 10 & 10 \\  
10 & 10 & 10  \\
9  & 10 & 10
\end{pmatrix}\, .
\end{equation}
The corresponding frequencies now are
$f_{1}=(\tfrac{10}{29},\tfrac{271}{841},\tfrac{280}{841})$ and
$f_{2}=(\tfrac{10}{29},\tfrac{1}{3},\tfrac{28}{87})$, and the growth
rates for these examples are at least $2^{1/29}$.

Our next examples use the generating words 
\begin{equation}
\begin{split}
w_{1}&\;=\;abcacbacabcbabcabacbcabcbacbcacba\qquad\hphantom{a}
\mbox{($n_{a}=11$, $n_{b}=11$, $n_{c}=11$)}\,,\\\
w_{2}&\;=\;abcacbcabacabcacbabcbacabacbcacba\qquad\hphantom{b}
\mbox{($n_{a}=12$, $n_{b}=10$, $n_{c}=11$)}\,,
\end{split}
\end{equation}
 of a Brinkhuis triple with $m=33$
\cite{G01}. Choosing as above 
$w_{a}^{(1)}=w_{1}$,
$w_{a}^{(2)}=\overline{w}_{1}$, 
$w_{b}^{(1)}=\sigma(w_2)$,
$w_{b}^{(2)}=\sigma(\overline{w}_2)$,
$w_{c}^{(1)}=\sigma^{2}(w_{\alpha})$ and
$w_{c}^{(2)}=\sigma^{2}(\overline{w}_{\alpha})$, where $\alpha\in\{1,2\}$, 
we obtain two Brinkhuis triples, this time with substitution matrices
$M_{\alpha}$ given by
\begin{equation}
M_{1}\;=\;
\begin{pmatrix}
11 & 11 & 11 \\  
11 & 12 & 11  \\
11 & 10 & 11
\end{pmatrix}\, ,\qquad
M_{2}\;=\;
\begin{pmatrix}
11 & 11 & 10 \\  
11 & 12 & 11  \\
11 & 10 & 12
\end{pmatrix}\, .
\end{equation}
The corresponding frequencies now are
$f_{1}=(\tfrac{1}{3},\tfrac{11}{32},\tfrac{31}{96})$
and $f_{2}=(\tfrac{331}{1024},\tfrac{11}{32},\tfrac{341}{1024})$.
Here, the growth rate is at least $2^{1/32}$.

Finally, we give one example with a rather large deviation from
equidistribution of letters.  This uses three generating words
\begin{equation}
\begin{split}
w_{1}&\;=\;abcacbacabacbcabacabcacbcabacbcacba\qquad\hphantom{bb}\!\!
\mbox{($n_{a}=13$, $n_{b}=10$, $n_{c}=12$)}\,,\\
w_{2}&\;=\;abcacbcabacbabcbacabcbabcabacbcacba\qquad\hphantom{ac}\!\!
\mbox{($n_{a}=12$, $n_{b}=12$, $n_{c}=11$)}\,,\\
w_{3}&\;=\;abcacbacabacbcabacabcbabcabacbcacba\qquad\hphantom{bc}\!\!
\mbox{($n_{a}=13$, $n_{b}=11$, $n_{c}=11$)}\,,
\end{split}
\end{equation}
of a Brinkhuis triple with $m=35$ \cite{G01}. Choosing 
$w_{a}^{(1)}=w_{1}$, 
$w_{a}^{(2)}=\overline{w}_{1}$,
$w_{b}^{(1)}=\sigma(w_2)$, 
$w_{b}^{(2)}=\sigma(\overline{w}_2)$,
$w_{c}^{(1)}=\sigma^{2}(w_{3})$ and 
$w_{c}^{(2)}=\sigma^{2}(\overline{w}_{3})$,
we obtain a Brinkhuis triple with substitution matrix
\begin{equation}
M\;=\;
\begin{pmatrix}
13 & 11 & 11 \\  
10 & 12 & 11  \\
12 & 12 & 13
\end{pmatrix}\,,
\end{equation}
which yields frequencies $f=(\tfrac{1}{3},\tfrac{16}{51},\tfrac{6}{17})$.
The growth rate is at least $2^{1/34}$.

To summarise, we proved the following.
\begin{lem}\label{th:freq}
The entropy of ternary square-free words with fixed letter 
frequency\/ $f_{a}$
is strictly positive for\/ $f_{a}\in
\{\tfrac{16}{51}, \tfrac{9}{28}, \tfrac{28}{87}, \tfrac{271}{841}, 
\tfrac{31}{96}, \tfrac{331}{1024}, \tfrac{280}{841}, \tfrac{341}{1024}, 
\tfrac{1}{3}, \tfrac{271}{812}, \tfrac{11}{32}, \tfrac{10}{29}, 
\tfrac{6}{17}\}$.
\qed 
\end{lem}

One should expect that the entropy is strictly positive for all letter
frequencies $f_{a}$ in an interval. However, it is not straightforward
to show that by using substitutions of Brinkhuis triples with
different letter frequencies. The reason is that, in general, the
infinite words obtained by such substitutions do not have well-defined
letter frequencies.

In the following sections, we are going to use methods from the theory
of generating functions and convex analysis \cite{St} which are often
applied in the context of statistical mechanics \cite{J00}. The free
energy of square-free words, which we will define below, is related to
the entropy function of square-free words with fixed letter density,
as follows from Proposition \ref{thm:ent}.  An immediate
consequence of the concavity of the entropy function is that the
entropy is strictly positive for all frequencies $f_{a}\in (16/51,
6/17)\approx (0.3137,0.3529)$, see below.

\section{Free energy}

Since the language of square-free words is subword closed, the
numbers\/ $s_{n,k}$ satisfy the submultiplicative inequality
\begin{equation}\label{form:submult}
s_{n+m,k} \;\le\; \sum_{l=0}^k\, s_{n,l}\, s_{m,k-l}\,.
\end{equation}
Consider the functions $s_n(q)$ defined by $s_n(q)\;=\;\sum_{k=0}^n\,
s_{n,k}\, q^k$.  These are polynomials in $q$ of degree not larger
than $n$.  The submultiplicative inequality (\ref{form:submult})
implies for the functions $s_n(q)$ that $s_{n+m}(q) \le s_n(q) \,
s_m(q)$ for $0<q<\infty$.  We are interested in the exponential growth
rate of $s_n(q)$.  To this end, define $F_n(q) := \frac{1}{n}\log
s_n(q)$.  The submultiplicative inequality yields
\cite[Lemma~A.1]{J00} that the limit $F(q):=\lim_{n \to \infty}
F_n(q)$ exists, and that $F(q) < \infty$ for $0<q<\infty$.  The
function $F(q)$ is called the {\it free energy}\/ of the model.  More
can be said about the properties of the free energy by using convexity
arguments.  These are largely independent of the underlying
combinatorial model and are discussed in detail in \cite[Sec.~2.1,
App.~B]{J00}.  We obtain
\begin{prop}
The functions\/ $F_n(q)=\frac{1}{n} \log s_n(q)$ of ternary
square-free words are continuous, analytic and convex in\/ $\log q$
in\/ $(0,\infty)$.  The free energy\/ $F(q)$ of ternary square-free
words
\begin{equation}
F(q)\;=\;\lim_{n \to \infty} F_{n}(q)
\end{equation}
exists and satisfies $F(q)<\infty$ for $q\in(0,\infty)$.  Moreover, it
is a convex function of\/ $\log q$ for $q\in(0,\infty)$.  If $F(q)$ is
finite, its right- and left-derivatives exist everywhere in\/
$(0,\infty)$, and they are non-decreasing functions of\/ $q$. The
function\/ $F(q)$ is differentiable almost everywhere, and wherever
the derivative\/ $dF(q)/dq$ exists, it is given by\/
$\lim_{n\to\infty}d F_n(q)/dq$.  \qed
\end{prop}

In the following, we will apply the results of the preceding section
in order to derive bounds on the free energy.  This will show that the
free energy $F(q)$ is finite for $0<q<\infty$.  Using the above
substitution rule (\ref{form:subs}) and the substitution rule given in
\cite{Tara02}, we first derive a lower bound on the free energy.
\begin{lem}
The free energy\/ $F(q)$ is bounded from below by
\begin{equation}
F(q) \;\ge\; 
\max \left\{ \frac{64}{233}\log q,\frac{13}{36}\log q \right\}\,.
\end{equation}
\end{lem}
\begin{proof}
Consider ternary square-free words $w_n$ of length $n=12k$, where
$k\in\mathbb N$, generated by the substitution rule (\ref{form:subs}),
with $w_1=c$.  Define $k_+(n)=13n/36 +\delta_+(n)$, which denotes the
number of letters of type $a$ in $w_n$.  Note that
$\delta_+(n)={o}(n)$.  We have $s_n(q)\ge s_{n,k_+(n)}q^{k_+(n)}$.
Taking the logarithm, dividing by $n$ and performing the limit leads
to $F(q) \ge \frac{13}{36} \log q$.  The second part of the statement
follows by the same argument with the substitution rule given in
\cite{Tara02}.
\end{proof}
\noindent {\bf Remark.} A weaker bound with $64/233$ replaced by
$11/36>64/233$ may be derived using the substitution
(\ref{form:subs}), where the role of $a$ and $b$ are interchanged.
\vspace{2ex}

We now turn to the question of an upper bound, which can be analysed
using the bounds for letter frequencies obtained in
\cite{Tara02,Tara03} or in Theorem \ref{th:fbounds}.
\begin{lem}
The free energy\/ $F(q)$ of ternary square-free words is bounded from
above by
\begin{equation}\label{form:ub}
F(q) \;\le\; -\log x_c \,+\, \max \left\{\frac{1780}{6481} \log q,
\frac{469}{1201} \log q \right\}
\end{equation}
where\/ $x_c=\lim_{n\to\infty} s_n^{1/n}\approx 0.768189$ is the
critical point of ternary square-free words.
\end{lem}
\begin{proof}
Assume that $q\neq1$. (The case $q=1$ has been discussed in Section 2,
where $F(1)=-\log x_c$ was proven.)  Assume that $B_n$ and $A_n$ are
numbers such that $s_{n,k}=0$ for $k>B_n$ or $k<A_n$, $s_{n,B_n}>0$,
and $s_{n,A_n}>0$.  For $1\neq q \in (0,\infty)$  we have the estimate
\begin{equation}
s_n(q)\;\le\; s_n\sum_{A_n}^{B_n}\,q^k 
\;=\; s_n\, \frac{q^{B_n+1}-q^{A_n}}{q-1}\,.
\end{equation}
Assume that $q>1$.  Taking the logarithm, dividing by $n$ and
performing the limit $n\to\infty$, this implies $F(q)\le \log x_c
+\epsilon_+ \log q$, where $\epsilon_+=\limsup_{n\to\infty} B_n/n$.
Note that $\epsilon_+\le 469/1201$, as follows from the bound given in
\cite{Tara03}.  A similar argument holds for $q<1$, involving the
lower bound $A_n$.  From \cite{Tara02}, we get the bound $1780/6481$.
Combining the two results, we get the inequality
(\ref{form:ub}).
\end{proof}
\noindent {\bf Remark.} A weaker bound with $(1780/6481, 469/1201)$
replaced by $(31/117, 39/97)$ follows from Theorem \ref{th:fbounds}.
\vspace{2ex}

Define the two-variable generating function $S(x,q)$
\begin{equation}
S(x,q)\; = \;\sum_{n=0}^\infty \sum_{k=0}^n\, s_{n,k}\, x^n\,  q^k 
\;=\;\sum_{n=0}^\infty\, s_n(q)\, x^n\,.
\end{equation}
Denote the radius of convergence of $S(x,q)$ by $x_c(q)$.  The curve
$x_c(q)$ is called {\it critical curve}, and the plot of $x_c(q)$ in
the $xq$-plane is called the {\it phase diagram}\/ of the model.  The
free energy is related to the critical curve by
\begin{equation}
x_c(q)^{-1} \;=\; \lim_{n\to\infty} s_n(q)^{1/n} \;=\; e^{F(q)}\,.
\end{equation}
We set $x_c=x_c(1)$ for the critical point of ternary square-free
words.  Bounds on the curve $x_c(q)$ can be derived from bounds on the
free energy $F(q)$ as given above.  This yields
\begin{equation}
x_{c}\,\min \{ q^{-1780/6481},q^{-469/1201}\} \;\le\; x_c(q) \;\le\; 
\min \{ q^{-64/233}, q^{-13/36}\}\,.
\end{equation}
The phase diagram is shown in Fig.~\ref{fig:phase}.  
\Bild{10}{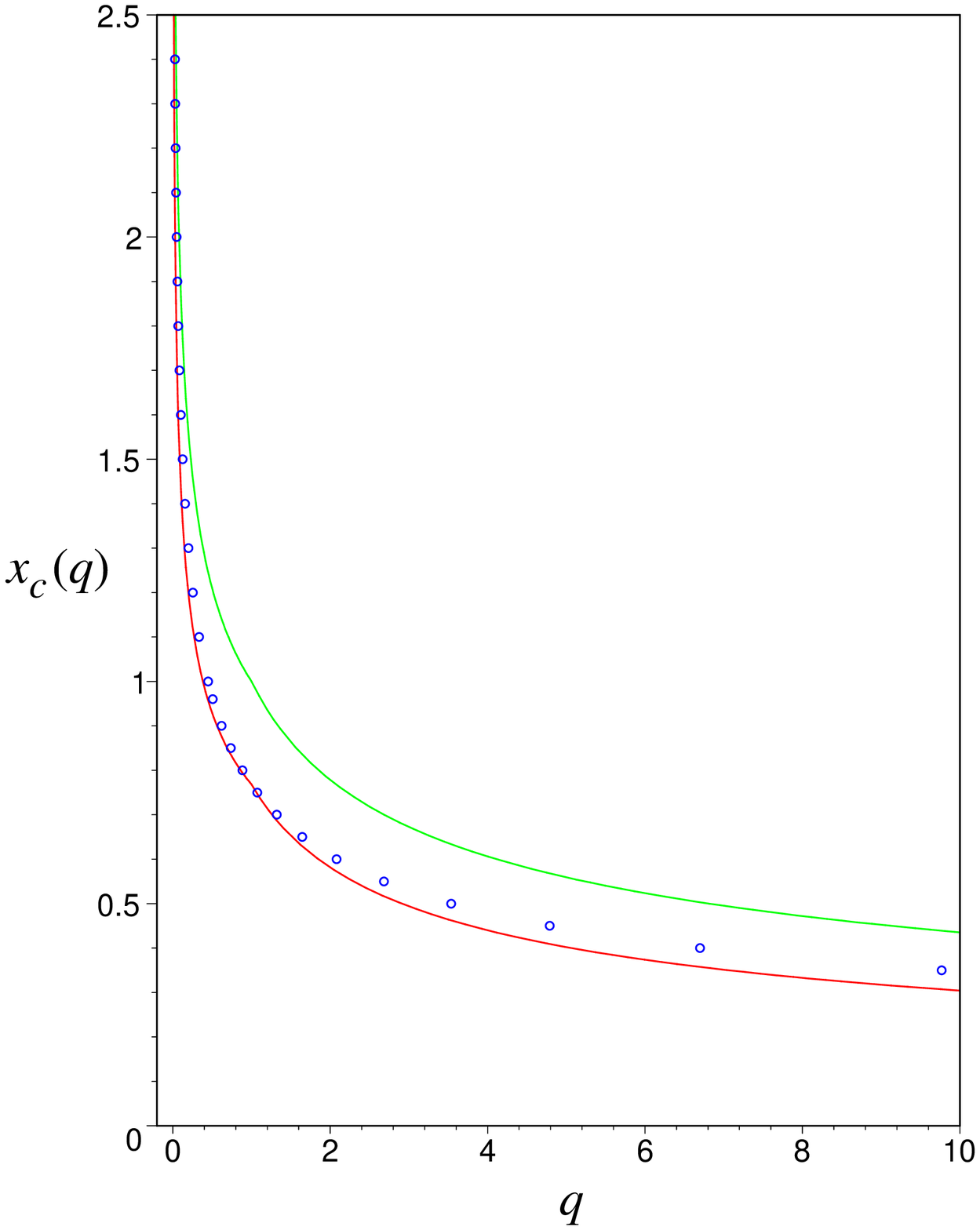}{Phase diagram of ternary square-free words, 
as extrapolated from exact enumeration data (circles). 
Upper and lower bounds on $x_c(q)$ are drawn for comparison.}{fig:phase}
Using the series data from exact enumeration for length $n\le 100$, we
extrapolated the values of $x_c(q)$ for different values of $q$, using
first order differential approximants \cite{G89}.  The critical curve
$x_c(q)$ is, within the analysed range of $q$, very close to the curve
$x_c \, q^{-1/3}$, reflecting the fact that the values $k=k(n)$ where
$s_{n,k}\neq0$ are sharply concentrated around $k=\lfloor n/3
\rfloor$.  For large values of $q$, such a form is, however, not
compatible with the derived bounds on $x_c(q)$.  Numerical analysis
suggests that the leading divergence of $S(x,q)$ is a simple pole,
which is approached uniformly in $x$ and $q$.  Thus, there is no
indication that the nature of the singularity changes, in contrast to
other examples from statistical mechanics, where such a change
indicates a phase transition \cite{J00}.

\section{Entropy and symmetry}

We now address the question of the number of ternary square-free
words, where we fix the frequency of letters of type $a$.  We consider
the number of square-free words $s_{n, \lfloor \epsilon n \rfloor}$ in
$n$ letters with $\lfloor \epsilon n \rfloor$ occurrences of the
letter $a$.  The number $\epsilon$ may thus be regarded as the
frequency of the letter $a$.  We are interested in the exponential
growth rate of $s_{n, \lfloor \epsilon n \rfloor}$.  This leads to the
question whether sequences of the form $\frac{1}{n} \log s_{n, \lfloor
\epsilon n \rfloor}$ have a limit as $n\to\infty$, which we then call
{\it entropy function}\/ $P(\epsilon)$. It is related to the free
energy $F(q)$ by a Legendre-Fenchel transform, as we will now show.

Note that there is a constant $K>0$ such that $0\le s_{n,k} \le K^n$
for each value of $n$ and $k$.  This follows from the existence of the
entropy $s$ of ternary square-free words.  Note also that there exists
a finite constant $C>0$, and numbers $A_n$ and $B_n$ such that
$s_{n,A_n}>0$ and $s_{n,B_n}>0$, and $s_{n,k}\ge 0$, when $0\le A_n <
k < B_n \le C n$.  This follows from the substitution rule
(\ref{form:subs}).  Take $A_n$ and $B_n$ such that $s_{n,k}=0$ if
$k<A_n$ or $k>B_n$.  Define the numbers
\begin{equation}
\epsilon_+ \;=\; \limsup_{n\to\infty}\, \frac{B_n}{n}\,, \qquad
\epsilon_- \;=\; \liminf_{n\to\infty}\, \frac{A_n}{n}\,.
\end{equation}
{}From \cite{Tara02,Tara03} and the substitution rule
(\ref{form:subs}), we have $0.361\approx 13/36 \le \epsilon_+ \le
469/1201\approx 0.391$ and $0.274649\approx 1780/6481\le \epsilon_-\le
64/233\approx 0.274678$.  Thus, the assumptions in
\cite[Thm.~3.19]{J00} are satisfied, and we obtain
\begin{prop}\label{thm:ent}
The entropy function\/ $P(\epsilon)$ of ternary square-free words exists
in\/ $(\epsilon_-,\epsilon_+)$ and is defined by
\begin{equation}\label{form:dens}
P(\epsilon)\;=\;\inf_{0<q<\infty} \{ F(q)-\epsilon \log q\}\,.
\end{equation}
Moreover, there is a sequence of integers\/ $\{ \sigma_n\}_{n=0}^\infty$
such that\/ $\sigma_n=o(n)$ and the limit
\begin{equation}
P(\epsilon)\;=\;\lim_{n\to\infty} \frac{1}{n} \log s_{n, \lfloor
\epsilon n \rfloor+\sigma_n}
\end{equation}
exists and is finite and concave in\/ $(\epsilon_-,\epsilon_+)$.
Lastly, note also that \/$\delta_n=\lfloor\epsilon n \rfloor+\sigma_n$
is the least value of\/ $k$ that maximises\/ $s_{n,k}\, {\tilde q}^k$,
where\/ $\tilde q$ is that value of\/ $q$ where the infimum is taken in\/
$(\ref{form:dens})$.  \qed
\end{prop}
\noindent {\bf Remark.} Together with Lemma 2, an immediate
consequence of the concavity of the entropy function is that the
entropy is strictly positive for all frequencies $\epsilon\in (16/51,
6/17)\approx (0.3137,0.3529)$.
\vspace{2ex}

We consider now the question where the entropy function takes its
maximum.  To this end, we assume a special regularity condition on the
free energy, whose validity is supported by the numerical analysis of
the preceding section, see also the discussion in the conclusion.
\begin{lem}\label{thm:ent2}
Let\/ $\epsilon\in (\epsilon_-,\epsilon_+)$.
If\/  $F(q)\in C^2(0,\infty)$, and if\/ $F(q)$ is strictly
convex in\/ $\log q$, we have\/ $P(\epsilon)\in
C^2(\epsilon_-,\epsilon_+)$ for the entropy function, and it is given by
\begin{equation}
P(\epsilon) \;=\; F\big(q(\epsilon)\big)\, -\, \epsilon\, 
\log q(\epsilon)\,,
\end{equation}
where\/ $q(\epsilon)$ is the unique positive solution of
\begin{equation}
\epsilon\;=\; q\,\frac{d}{dq}F(q)\,.
\end{equation}
The entropy function\/ $P(\epsilon)$ attains its global maximum at\/ $q=1$.
\end{lem}
\begin{proof}
Since $F(q)$ is convex in $\log q$ and continuous, and $F(q)\ge\max\{
\epsilon_-\log q, \epsilon_+\log q\}$, the infimum in
(\ref{form:dens}) occurs at a unique value $q=q(\epsilon)\in
(0,\infty)$.  Since $F(q)\in C^1(0,\infty)$, we obtain $\epsilon=q
F'(q)=\frac{d}{d(\log q)}F(q)$ as an implicit equation for
$q(\epsilon)$.  This uniquely defines a positive function
$q=q(\epsilon)\in C^1(\epsilon_-,\epsilon_+)$, since strict convexity
of $F(q)$ and $F(q)\in C^2(0,\infty)$ implies $\frac{d^2}{d(\log
q)^2}F(q)\neq0$.  We have explicitly $P'(\epsilon)=-\log q(\epsilon)$,
which shows that $P(\epsilon)\in C^2(\epsilon_-,\epsilon_+)$, and
$-\infty<P''(\epsilon)=-(\frac{d^2}{d(\log q)^2}F(q))^{-1}<0$.  This
implies that $q=1$ is a local maximum of $P(\epsilon)$.  Due to the
concavity of $P(\epsilon)$, it is the global maximum.
\end{proof}

We note that at $q=1$, the letter density $\epsilon=F'(1)$ is the mean
letter density, which was determined above to be $\epsilon=1/3$ by a
symmetry argument.  Thus, under the above regularity assumption,
maximum entropy occurs at equal (mean) letter density
$\epsilon_a=\epsilon_b=\epsilon_c=1/3$.  This is an example of the
more general result that maximum entropy in occurs at points of
maximum symmetry, see \cite{RHHB98} for the concept of symmetry and
its implications for the free energy and entropy of the combinatorial
problem of random tilings, which is applicable in this case.

\section{Conclusions}

In this article, we considered the growth rate, or the entropy, of the
set of ternary square-free words. By computing generating functions
$S^{(\ell)}(x)$ for length-$\ell$ square-free words, where the
condition of square-freeness is truncated at length $\ell$, we
verified an upper bound proposed in \cite{NZ99} and slightly improved
it. The pattern of poles of these generating functions, and their
behaviour as $\ell$ increases, points towards a natural boundary for
the generating function $S(x)$.

The presence of a natural boundary in a model indicates that it cannot
be solved exactly in terms of standard functions of mathematical
physics, which obey linear differential equations with polynomial
coefficients \cite{G00}.  This would exclude, for ternary square-free
words, an exact value for the entropy and the functional form of the
free energy.  It may even be difficult to prove the existence of a
critical exponent, compare the related self-avoiding walk problem
\cite{J00}.

In the ternary alphabet, no letter is preferred by the condition of
square-freeness. Thus, averaging over the entire sets of ternary
square-free words, all letters appear equally often. However, in a
single infinite word this need not be the case, indeed, the letter
frequency may not be well-defined. However, one can derive limits on
the minimum or maximum frequency of a given letter in an infinite
ternary square-free words, and by explicitly constructing infinite
words with given well-defined frequencies by means of substitution
rules the minimum and maximum frequency can be bounded from above and
below. We obtained limits from counting square-free words up to a
certain length, sharper limits were given recently in
\cite{Tara02,Tara03}.  The bounds for the maximum frequency can
certainly be further improved employing the approach of
\cite{KK97,KKT98,Tara02}.

Lower bounds on the entropy are based on Brinkhuis triples and their
generalisations. We used these to prove that, for a list of rational
values, the entropy of the set of square-free words with a fixed
letter frequency is strictly positive.  Together with the concavity of
the entropy function, obtained by methods of convex analysis and
statistical mechanics, this led to the result that the entropy is
strictly positive on an entire interval.

Concerning the entropy function, it would be interesting to extend the
interval of strict positivity by providing sharper bounds from
suitable substitution rules. This might be achievable by following and
suitably modifying the approach taken in \cite{KK97,KKT98,Tara02}.  It
is conceivable, albeit not necessary, that there exists a region of
frequencies for which infinite square-free words exist, but the
entropy vanishes, because the number of square-free words with that
given letter frequency grows sub-exponentially. Such behaviour has
been reported for $k$\/th-power-free binary square-free words with
rational powers in the range $2<k\le$7/3 \cite{G01b}.

Further, it is necessary to prove the validity of the regularity
assumption on the free energy in Theorem \ref{thm:ent2}. In contrast
to other problems in statistical mechanics \cite{J00}, there is no
indication of a phase transition in the model of ternary square-free
words, wherefore an analytic free energy is expected.

It would also be interesting to analyse the letter distribution using
probabilistic methods.  Similar examples lead, in an appropriate
scaling limit, to Gaussian distribution functions \cite{L96}.

\section{Acknowledgements}

We thank the Erwin Schr\"{o}dinger International Institute for
Mathematical Physics in Vienna for support during a stay in winter
2002/2003, where part of this work was done. CR would like to
acknowledge financial support by the German Research Council (DFG).
We are grateful to Jeffrey O. Shallit and to Bernd Sing for making us
aware of references \cite{Tara02} and \cite{Z58}, respectively.
Furthermore, we thank R.~Kolpakov and Y.~Tarannikov for providing us
with their recent (partly unpublished) results. UG also wishes to
acknowledge useful discussions with X.~Y.~Sun and D.\ Zeilberger.

\end{document}